\documentclass[onecolumn,10pt]{IEEEtran}

\usepackage{graphicx}

\usepackage{hyperref}%
\hypersetup{colorlinks=true, linkcolor=blue, breaklinks=true, urlcolor=blue}
\usepackage{url,doi}

\usepackage[numbers]{natbib}
\usepackage{url}
\usepackage{bbm}
\usepackage{bbold}
\usepackage{bm}
\usepackage{amsmath}
\usepackage{xcolor}

\usepackage{subfigure}

\usepackage{tikz}
\usetikzlibrary{automata,positioning}
\usetikzlibrary{calc}
\usetikzlibrary{arrows.meta}

\newcommand{\oneton}{\{1,\dots,n\}}
\newcommand{\vecz}{\textup{vec}}
\newcommand{\diag}{\textup{diag}}

\newcommand{\subscr}[2]{#1_{\textup{#2}}}
\newcommand{\Hrt}{\subscr{\mathcal{H}}{ret-time}}
\newcommand{\Hrate}{\subscr{\mathcal{H}}{rate}}

\newtheorem{problem}{Problem}

\title{Markov Chain-Based Stochastic Strategies \\ for Robotic Surveillance}

\author{Xiaoming Duan and Francesco Bullo
	\thanks{Department of Mechanical
		Engineering and the Center of Control, Dynamical
		Systems and Computation, UC Santa Barbara, CA 93106-5070, USA. email: xiaomingduan@ucsb.edu, bullo@ucsb.edu.}
}
\begin{document}


\maketitle

\begin{abstract}
This article surveys recent advancements of strategy designs for persistent robotic surveillance tasks with the focus on stochastic approaches. The problem describes how mobile robots stochastically patrol a graph in an efficient way where the efficiency is defined with respect to relevant underlying performance metrics. We first start by reviewing the basics of Markov chains, which is the primary motion model for stochastic robotic surveillance. Then two main criteria regarding the speed and unpredictability of surveillance strategies are discussed. The central objects that appear throughout the treatment is the hitting times of Markov chains, their distributions and expectations. We formulate various optimization problems based on the concerned metrics in different scenarios and establish their respective properties. 
\end{abstract}



\section{INTRODUCTION}
This survey presents a stochastic approach for the design of robotic
surveillance algorithms. We adopt Markov chains as the main algorithmic
building block and discuss various properties relevant in surveillance
applications. This approach leads to easily implementable and intrinsically
unpredictable surveillance algorithms in which robots' visit frequency at
different locations can be conveniently specified in a probabilistic
fashion. In particular, the unpredictability of a random walk over a graph
is a highly desirable characteristic for surveillance strategies in
adversarial settings, because it makes it hard for potential intruders to
take advantage of surveillance robots' motion patterns.

In persistent surveillance tasks, mobile robots patrol and visit sites of
interest in the environment to collect sensing data, detect anomalies or
intercept intrusions. Persistent surveillance is a key component in
numerous applications such as environmental
monitoring~\cite{CY-YZ-ZL:15,AS-HK:16}, infrastructure
inspection~\cite{DDP-AM-GC-AD:10,SW-JCC-JM-JC-FSM-PUL-MV:17} and urban
security~\cite{TDR:10}. Different from traditional coverage or search and
exploration problems, surveillance tasks require robots to repeatedly visit
different locations and achieve continuous monitoring. Surveillance
strategies specify in which order and how often the mobile robots should
visit different locations so that desirable performance is maintained. In
this survey, we adopt a model where the physical environment is discretized and
modeled by a connected graph: each node of the graph represents a location
of interest in the environment and the edges in the graph represent the
paths between locations. Scenarios where mobile robots move in continuous
space are studied for example in
\cite{FP-JWD-FB:11h,CGC-XL-XCD:13,NZ-XY-SBA-CGC:18,YW-YW-XL-NZ-CGC:19}. There
are also works that consider (stochastic) dynamics at the locations of
interest such
as~\cite{SLS-MS-DR:11,JY-SK-DR:15,NR-SSK:19,XY-SBA-NZ-CGC:20}.

\paragraph*{Deterministic strategy}
In the design of deterministic surveillance strategies, the performance
criterion of typical interest is \emph{idleness}, also known as refresh
time or latency. At any given time instant, the idleness of a location is
the time elapsed since the last time the location was visited by any
surveillance robots. The \emph{worst idleness} (respectively, \emph{average
  idleness}) at that time instant is the largest (respectively, average)
value of idleness for all locations. Then, one seeks a strategy that
minimizes the largest worst (average) idleness over the entire (possibly
infinite) surveillance period \cite{AM-GR-JDZ-AD:03}. In the early
work~\cite{YC:04}, the author shows that the cyclic strategy where a robot
travels with the maximum speed on a shortest closed route through all
locations is optimal in the sense that it minimizes the maximum worst
idleness. The author also proposes a partition-based approach for the
multi-robot case where the environment is partitioned into disjoint regions
so that each robot patrols a single region. In~\cite{YE-NA-GAK:09}, a
multi-robot system composed of homogeneous robots is placed on a cyclic
route in the environment, and a uniform visit frequency to all locations is
achieved by carefully designing the displacements between robots. The
authors of \cite{fp-af-fb:09v} study team strategies with minimum refresh
time for the line, tree and cyclic graphs over a finite surveillance
period; optimal and constant factor suboptimal strategies are
obtained. When locations have different priorities, one can define the
natural notion of \emph{weighted worst idleness}, whereby idleness is
penalized proportionally to the location priority. Single-robot
surveillance strategies that minimize the maximum weighted idleness are
studied in \cite{SA-EF-SLS:14}. The authors design two approximate
algorithms to compute strategies whose costs are within factors of the
optimal cost; these factors depend on the distribution of weights or the
dimension of the graph. Given constraints on the weighted idleness at all
locations, the authors in \cite{ABA-SLS:19} propose an approximate
algorithm to compute the minimum number of robots required to satisfy the
idleness requirements. The recent
work~\cite{RA-MDB-KB-JG-ML-AN-BR-RS-HW-HY:20} proposes surveillance plans
for a fixed number of robots that minimize the maximum weighted
idleness. In summary, optimization problems concerning the idleness metrics
and deterministic surveillance strategies are usually of combinatorial
nature and, in most cases, approximation algorithms and suboptimal
solutions are sought.

While clearly leaving any location of interest in the environment
unattended for an extended period of time should be avoided, there are two
main challenges for deterministic surveillance strategies: (i) as pointed
out in \cite{KS-DMS-MWS:09}, it is not always possible to assign an
arbitrary visit frequency to different locations in the graph (see recent
related work in~\cite{DK-LJS-OT:18}); (ii) deterministic strategies are
fully predictable and can be easily exploited by potential intruders in
adversarial settings.  In this case, stochastic strategies become
particularly appealing.

\paragraph*{Stochastic strategy}
In adversarial settings, potential intruders may be able to observe and
learn the strategies of the surveillance robot and devise intrusion plans
accordingly. For example, when surveillance robot adopts a deterministic
strategy, the intruder can confidently attack a location immediately after
the surveillance robot leaves that location, because it knows for certain
that the robot will not return to that location for a deterministic known
duration of time. In such scenarios, robotic surveillance problems have to
be tackled by resorting to randomized approaches. One popular technique in
such approaches is to describe the surveillance strategy as a Markov
chain. There are several advantages to using Markov chains as stochastic
strategies, including: 1) the intrinsic stochasticity of Markov chains
leads to more effective strategies against rational intruders; 2) the visit
frequency to different locations can be easily assigned through the
stationary distribution; and 3) Markov chains as routing algorithms are
lightweight and require minimal amount of resources to implement and
execute. Although high-order Markov chains are more versatile, their state
space grows exponentially fast with the order (memory length), which
results in high computational complexity or intractability. Therefore,
first-order Markov chains are the default choices in the design. Depending
on whether an intruder model is explicitly specified or not, there are two
different common formulations in the literature. When no intruder model is
assumed, stochastic strategies are designed based on performance metrics
such as \emph{speed} and \emph{unpredictability}; when a specific intruder
model is adopted, tailored strategies to the intruder behaviors are
carefully constructed.

\emph{Metric-based design:} Commonly-used design metrics for stochastic
strategies in robotic surveillance include visit frequency, speed and
unpredictability. In \cite{GC-AS:11}, the authors propose a Markov
chain-based algorithm for a team of robots to achieve continuous coverage
of an environment with a desirable visit frequency. Smart nodes deployed at
different locations are responsible for recording the visits by robots and
directing robots to the next site to be visited. To obtain fast strategies
for the surveillance agent, the authors of \cite{RP-PA-FB:14b} study Markov
chains with the minimum weighted mean hitting time. Here, travel times
between different locations are given as edge weights on the graph. The
minimization of the weighted mean hitting time is transcribed into a convex
optimization problem for the special case where candidate Markov chains are
assumed reversible. The problem also has a semidefinite reformulation and
can be solved efficiently. A similar notion of mean hitting time for a
multi-robot system is proposed and studied in \cite{RP-AC-FB:14k}. To
obtain maximally unpredictable surveillance strategies, the authors in
\cite{MG-SJ-FB:17b} design Markov chains with maximum entropy rate:
numerous properties of the maxentropic Markov chain are established and a
fast algorithm to compute the optimal strategy is proposed.
A new notion that quantifies unpredictability of surveillance strategies
through the entropy in return time distributions has been recently
introduced and characterized in~\cite{XD-MG-FB:17o}. This new concept is
particularly relevant and useful in cases when only local observations are
available to potential intruders. Similar studies on introducing randomness
in return times also appeared in \cite{CDA-NB-SC:19,NB-SC:20}, where the
authors propose to insert random delays into the return times.

\emph{Intruder model-based design:} In contrast to the metric-based designs
where no explicit intruder behaviors are assumed, intruder model-based
designs leverage available information on the intruder to achieve
improved/guaranteed performance. In~\cite{NA-GK-SK:11}, a multi-robot
perimeter patrol problem is proposed where a randomized strategy is used to
defend against a strategic adversary. The adversary knows the surveillance
strategy as well as the robots' locations on the perimeter, and it aims to
maximize its probability of success by attacking the weakest spot along the
perimeter.
Coordinated sequential attacks in the multi-robot perimeter patrol problem
are considered in \cite{ESL-NA-SK:19}. The authors in \cite{NB-NG-FA:12}
define a patrolling security game where a Markovian surveillance robot
tries to capture an intruder that knows perfectly the location and the
strategy of the surveillance agent. The intruder has the freedom to choose
where and when to attack so that the probability of being captured is
minimized. Various intruder models that characterize different manners in
which the intruders attack the locations are discussed
in~\cite{ABA-SLS:16}. The authors design pattern search algorithms to solve
for a Markov chain that minimizes the expected reward for the
intruders. Variations of intruder models where intruders have limited
observation time and varying attack duration are studied in
\cite{ABA-SLS:18} and \cite{HY-ST-KSL-SL-JG:19}, respectively. Many of the
aforementioned models can be formalized as Stackelberg games. The
Stackelberg game framework has been successfully applied in practical
security domain applications~\cite{AS-FF-BA-CK-MT:18}. However, only
limited progress on efficient computational methods with optimality
guarantees has been made.

In this paper, we focus on the metric-based design of stochastic
surveillance strategies. Concepts and notions discussed here also appear in
other settings in the literature.

\paragraph*{Organization}  We first review the basics of Markov chains and relevant design metrics in Section~\ref{sec:MCs}. Then, fast and unpredictable surveillance strategies in various scenarios are discussed in Section~\ref{sec:fastMC} and Section~\ref{sec:unpredictableMC}, respectively. We conclude the paper in Section~\ref{sec:conclusion}.

\paragraph*{Notation} Let $\mathbb{R}$, $\mathbb{R}^n$ and $\mathbb{R}^{m\times n}$ be the set of real numbers, real vectors of dimension $n$ and real matrices, respectively. The set of elementwise positive vectors and nonegative matrices are denoted by $\mathbb{R}^n_{>0}$ and $\mathbb{R}^{m\times n}_{\geq0}$. We denote the vector of $1$'s in $\mathbb{R}^n$ by $\mathbb{1}_n$ and the $i$-th standard unit vector by $\mathbb{e}_i$, whose dimension will be clear from the context. For a vector $v\in\mathbb{R}^n$, $\diag(v)\in\mathbb{R}^{n\times n}$ is a diagonal matrix with diagonal elements being $v$; for a matrix $S\in\mathbb{R}^{n\times n}$, $\diag(S)\in\mathbb{R}^{n\times n}$ is a diagonal matrix with diagonal elements being the same as that of $S$. For a matrix $S\in\mathbb{R}^{m\times n}$, the vectorization $\vecz(S)\in\mathbb{R}^{mn}$ is constructed by stacking the columns of $S$ on top of one another. The indicator function is denoted by $\mathbf{1}_{\{\cdot\}}$.

\section{MARKOV CHAINS AND RELEVANT METRICS}\label{sec:MCs}
In this section, we review the basics of discrete-time homogeneous Markov chains and the relevant metrics in robotic surveillance applications.
\subsection{Markov Chains}
A first-order discrete-time homogeneous Markov chain defined over the state space $\oneton$ is a sequence of random variables $X_k$ for $k\geq0$ that satisfies the Markov property: $\mathbb{P}(X_{k+1}=i_{k+1}\,|\,X_{k}=i_k,\dots,X_0=i_0)=\mathbb{P}(X_{k+1}=i_{k+1}\,|\,X_{k}=i_k)$ for all $k\geq0$ and $i_0,\dots,i_{k+1}\in\oneton$. The Markov chain $\{X_k\}_{k\geq0}$ has an associated row-stochastic transition matrix $P\in\mathbb{R}^{n\times n}$ such that the $(i,j)$-th element  $p_{ij}=\mathbb{P}(X_{k+1}=j\,|\,X_{k}=i)$ for $i,j\in\oneton$. The transition diagram of a Markov chain is a directed graph $\mathcal{G}=(V,\mathcal{E})$, where $V=\oneton$ is the set of nodes and $\mathcal{E}\subset V\times V$ is the set of edges, and $(i,j)\in\mathcal{E}$ if and only if $p_{ij}>0$. A Markov chain is irreducible if its transition diagram is strongly connected. A finite-state irreducible Markov chain $P$ has a unique stationary distribution $\bm{\pi}\in\mathbb{R}^n$ that satisfies $\bm{\pi}^\top P=\bm{\pi}^\top$ and $\bm{\pi}_i>0$ for $i\in\{1,\dots,n\}$. Moreover, the stationary distribution $\bm{\pi}$ has the interpretation that regardless of the initial condition \cite[Theorem 2.1]{DA-JAF:02}, 
\begin{equation}\label{eq:sd}
\frac{1}{t+1}\sum_{k=0}^{t}\mathbf{1}_{\{X_k=i\}}\xrightarrow[]{\textup{as } t\rightarrow\infty}\bm{\pi}_i\quad\textup{almost surely.}
\end{equation}
The physical meaning of~\eqref{eq:sd} is that for an irreducible Markov chain, the stationary distribution encodes the visit frequency to different states in the long run. A Markov chain is reversible if $\bm{\pi}_ip_{ij}=\bm{\pi}_jp_{ji}$ for all $i,j\in\{1,\dots,n\}$. One can verify that if a positive vector $\bm{\pi}\in\mathbb{R}^n_{>0}$ satisfies the reversibility condition for a transition matrix $P$, then it must be a stationary distribution of $P$. We refer the readers to \cite{JGK-JLS:76,JRN:97} for more about Markov chains. We focus on discrete-time homogeneous irreducible Markov chains with finite states in this article.

In robotic surveillance applications, the environment is sometimes modeled by a weighted directed graph $\mathcal{G}=(V,\mathcal{E},W)$ where the elements of the weight matrix $W=[w_{ij}]$ represent the travel distances or times between locations. A Markov chain defined over a weighted directed graph $\mathcal{G}=(V,\mathcal{E},W)$ has the transition diagram that is consistent with $\mathcal{G}$, i.e., the transition probability $p_{ij}\geq0$ if $(i,j)\in\mathcal{E}$ and $p_{ij}=0$ otherwise, for all $i,j\in\{
1,\dots,n\}$. Moreover, in this case, the weights specify the amounts of time it takes for the Markov chain to transition from one state to another.

\subsection{Hitting Times}
For a finite-state discrete-time Markov chain, the first hitting time from state $i$ to state $j$ is a random variable defined by
\begin{equation}\label{eq:hittingtime}
T_{ij}=\min\{k\,|\,X_0=i,X_k=j,k\geq1\}.
\end{equation}
This notion can be extended to Markov chains defined over a weighted graph $\mathcal{G}=(V,\mathcal{E},W)$, in which case the hitting time $T_{ij}^{\textup{w}}$ satisfies
\begin{equation}\label{eq:hittingtimeweight}
T_{ij}^{\textup{w}}=\min\{\sum_{\ell=0}^{k-1}w_{X_{\ell},X_{\ell+1}}\,|\,X_0=i,X_k=j,k\geq1\}.
\end{equation}
The hitting times measure how fast a Markov chain transitions from one state to another. If we require the Markov chains to be as fast as possible, then the ``speed'' can be conveniently quantified by hitting times and their expectations.

\subsection{Entropy}
Although Markov chains usually come with unpredictability, different Markov chains have different level of randomness as measured by some notion of entropy. For example, a Hamiltonian tour in a graph (if there exists one) can be represented by a Markov chain with transition matrix being an irreducible permutation matrix. However, a Hamiltonian tour is clearly deterministic and thus fully predictable. We choose to use entropy to quantify unpredictability of Markov chains as: i) entropy is a well-established fundamental concept that
characterizes randomness of probability distributions; ii) when the surveillance agent is highly entropic, its strategy (motion pattern) is potentially hard for the intruders to learn.

The entropy $\mathcal{H}(X)$ of a discrete random variable $X$ taking values in $\{1,\dots,n\}$ with distribution $\mathbb{P}(X=i)=p_i$ for $i\in\oneton$ is defined by
\begin{equation*}
\mathcal{H}(X)=-\sum_{i=1}^np_i\log p_i.
\end{equation*}
In order to quantify the randomness in Markovian strategies, two different notions of entropy have been used in the literature.
\subsubsection{Entropy Rate}\label{sec:entropyrateintro}
A classic notion that measures the randomness of a Markov chain is the entropy rate. For a Markov chain with the initial distribution being its stationary distribution, its entropy rate is defined by \cite[Theorem 4.2.4]{TMC-JAT:12}
\begin{equation}\label{eq:entropyrate}
\Hrate(P)=-\sum_{i=1}^n\bm{\pi}_i\sum_{j=1}^{n}p_{ij}\log p_{ij}.
\end{equation}
The entropy rate measures the randomness in the sequence of states visited by a Markov chain \cite{MG-SJ-FB:17b, XD-MG-FB:17o}.

\subsubsection{Return Time Entropy}
The entropy of Markov trajectories is first proposed and studied in~\cite{LE-TMC:93}, where a trajectory $\mathcal{T}_{ij}$ from a state $i$ to a state $j$ is a path with initial state $i$, final state $j$ and no other states being $j$. It is shown that the entropy of the random return trajectory $\mathcal{H}(\mathcal{T}_{ii})=\bm{\pi}_i\Hrate(P)$ for all $i\in\{1,\dots,n\}$. In robotic surveillance applications, we are more concerned with the amount of time it takes for the surveillance agent to come back instead of the exact trajectory followed by the surveillance agent. This is because often times, intruders may have access to the local observations of the inter-visit times to a location, and they could attack accordingly based on this information. We thereby define the return time entropy of a Markov chain as a weighted sum of the entropy of return times to locations with weights being the stationary distribution, i.e.,
\begin{equation}\label{eq:returntimeentropy}
\Hrt(P)=\sum_{i=1}^{n}\bm{\pi}_i\mathcal{H}(T_{ii}),
\end{equation}
where $T_{ii}$ is the return time (a discrete random variable) defined by~\eqref{eq:hittingtime}. We will see how this definition can be generalized to the cases when there are travel times on the surveillance graph.

\subsection{Running Examples}
\begin{figure}[http]\centering
  \subfigure[$3\times3$ grid graph]{\includegraphics[width=0.35\textwidth]{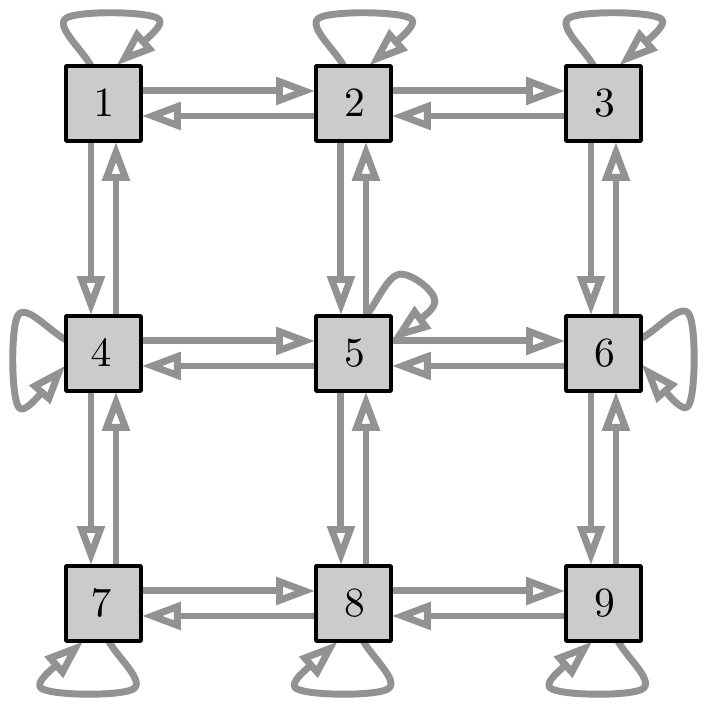}} \hfil
  \subfigure[Map of San Francisco]{\includegraphics[width=0.3\textwidth]{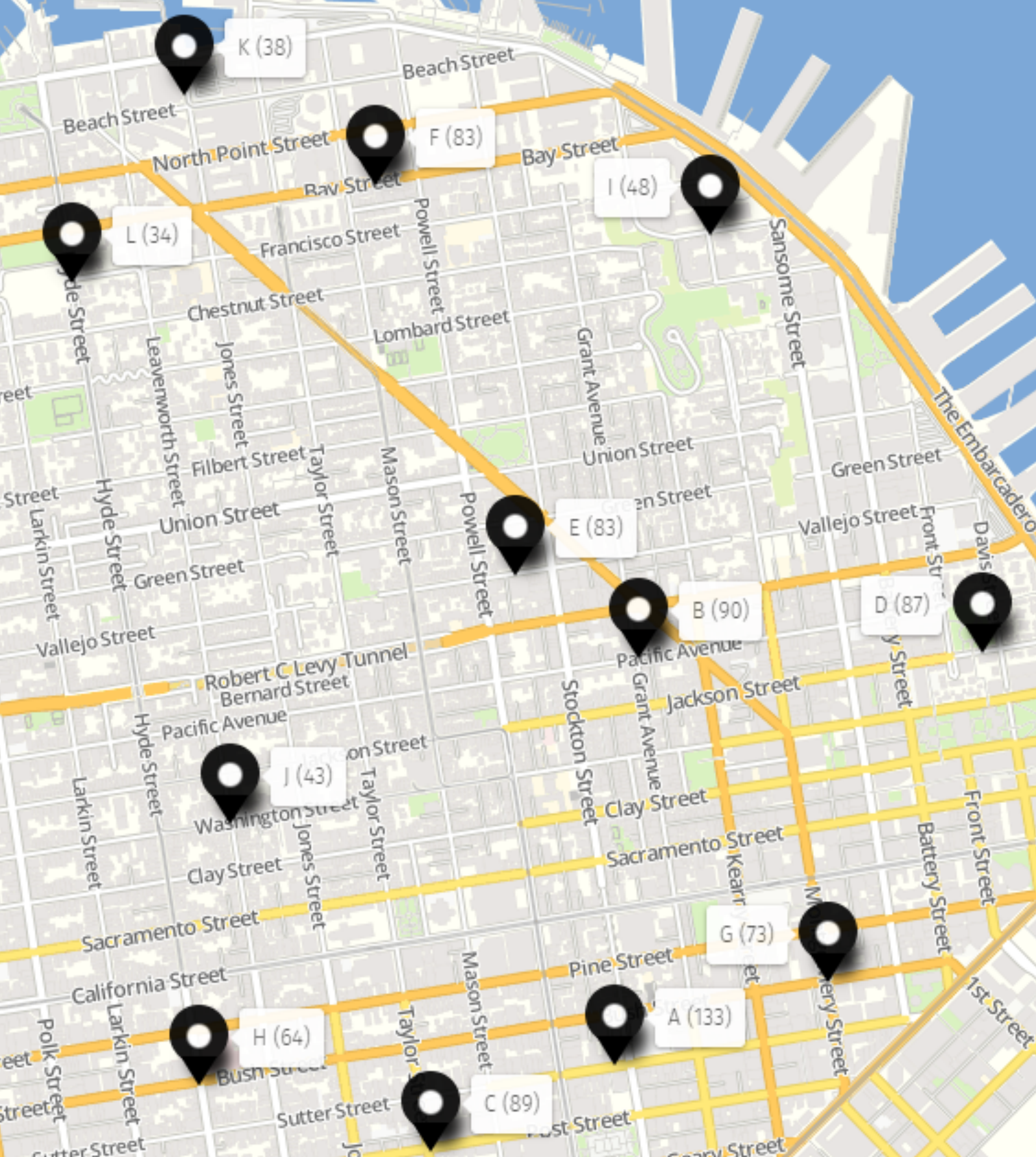}}
  \caption{Two running examples: (a)  a $3\times3$ grid graph with self loops at all locations; (b) a region map of
    San Francisco (created with uMap \cite{uMap:13}) modeled by a fully
    connected graph where the number of crimes during a time period is
    recorded at $12$ locations.}\label{fig:runningexample}
\end{figure}

Before we dive into the details of the design of Markov strategies, we introduce two running examples with which we will illustrate and visualize the optimal strategies to be discussed in later sections. We consider a $3\times3$ grid graph and a region map of San Francisco (SF)~\cite[Section 6.2]{SA-EF-SLS:14} as shown in Fig~\ref{fig:runningexample}. In the grid graph, we assume unit travel times (unweighted) and the surveillance agent is expected to visit all locations with the same frequency, i.e., $\bm{\pi}=\frac{1}{9}\mathbb{1}_9$. In the SF map, there are $12$ locations forming a fully connected weighted graph. The weights on the edges, as documented in Table~\ref{tb:pairwiseTime}, are the quantized by-car travel times (in minutes) between locations. Moreover, each location is associated with a number that indicates the number of crimes recorded at that location during a certain period of time. The surveillance agent is expected to have a visit frequency to different locations proportional to the crime rates, i.e.,  $\bm{\pi}=\begin{bmatrix}
\frac{133}{866}, \frac{90}{866},\frac{89}{866},\frac{87}{866},
\frac{83}{866}, \frac{83}{866},\frac{74}{866}, \frac{64}{866},
\frac{48}{866}, \frac{43}{866}, \frac{38}{866}, \frac{34}{866}
\end{bmatrix}^\top$.

\begin{table}[http]
	\centering
	\caption{Quantized pairwise by-car travel times in minutes}\label{tb:pairwiseTime}
	\begin{tabular}{ccccccccccccc}
		\hline\rule{0pt}{1\normalbaselineskip}
		Location&A&B&C&D&E&F&G&H&I&J&K&L \\ \hline\rule{0pt}{1\normalbaselineskip}A&1&	3&	3&	5&	4&	6&	3&	5&	7&	4&	6&	6\\
		B&3&	1&	5&	4&	2&	4&	4&	5&	5&	3&	5&	5\\
		C&3&	5&	1&	7&	6&	8&	3&	4&	9&	4&	8&	7\\
		D&6&	4&	7&	1&	5&	6&	4&	7&	5&	6&	6&	7\\
		E&4&	3&	6&	5&	1&	3&	5&	5&	6&	3&	4&	4\\
		F&6&	4&	8&	5&	3&	1&	6&	7&	3&	6&	2&	3\\
		G&2&	5&	3&	5&	6&	7&	1&	5&	7&	5&	7&	8\\
		H&3&	5&	2&	7&	6&	7&	3&	1&	9&	3&	7&	5\\
		I&8&	6&	9&	4&	6&	4&	6&	9&	1&	8&	5&	7\\
		J&4&	3&	4&	6&	3&	5&	5&	3&	7&	1&	5&	3\\
		K&6&	4&	8&	6&	4&	2&	6&	6&	4&	5&	1&	3\\
		L&6&	4&	6&	6&	3&	3&	6&	4&	5&	3&	2&	1  \\\hline
	\end{tabular}
\end{table}

\section{FAST SURVEILLANCE: MINIMIZATION OF MEAN HITTING TIME}\label{sec:fastMC}
As hitting times measure how fast a Markov chain travels on a graph, it is natural to formulate a hitting time minimization problem for robotic surveillance scenarios where the surveillance agent is expected to visit different locations in the graph as fast as possible. We first consider the case when the travel times are unitary and introduce the well-known Kemeny's constant. Then, we investigate how the formulation changes when arbitrary positive real-valued travel times are present. The development in Subsections~\ref{subsec:hittingunit} and~\ref{subsec:hittingarb} follows from~\cite{RP-PA-FB:14b}; the multi-vehicle case in Subsection~\ref{subsec:multihitting} is studied in~\cite{RP-AC-FB:14k}; the meeting times of Markovian agents in Subsection~\ref{subsec:meetingtime} are discussed in~\cite{XD-MG-RP-FB:14l}.

\subsection{Unit Travel Times}\label{subsec:hittingunit}
Let $m_{ij}$ be the expectation of the hitting time $T_{ij}$ defined in~\eqref{eq:hittingtime} for an irreducible Markov chain $P$, i.e., $m_{ij}=\mathbb{E}[T_{ij}]$, then the Markov property implies the following recursive relation
\begin{equation}\label{eq:meanhitting}
m_{ij}=p_{ij} + \sum_{k\neq j}p_{ik}(1+m_{kj})=1+\sum_{k=1}^np_{ik}m_{kj}-p_{ij}m_{jj}.
\end{equation}
We can write \eqref{eq:meanhitting} in matrix form as
\begin{equation}\label{eq:meanhittingmatrix}
M=\mathbb{1}_n\mathbb{1}_n^\top+P(M-\diag(M)).
\end{equation}
Two observations are immediate from~\eqref{eq:meanhittingmatrix}:
\begin{enumerate}
\item\label{enu:diagonal} by multiplying from the left the stationary distribution $\bm{\pi}^\top$ of $P$ on both sides of~\eqref{eq:meanhittingmatrix} and using $\bm{\pi}^\top P=\bm{\pi}^\top$, we have $m_{ii}=\frac{1}{\bm{\pi}_i}$ for all $i\in\{1,\dots,n\}$ or $\diag(M)=\diag(\bm{\pi})^{-1}$;

\item by multiplying from the right $\bm{\pi}$ on both sides of~\eqref{eq:meanhittingmatrix} and using $\diag(M)=\diag(\bm{\pi})^{-1}$, we have $(I-P)M\bm{\pi}=0$, which implies that $M\bm{\pi}$ is a constant multiple of $\mathbb{1}_n$ given $P$ is irreducible. 
\end{enumerate}
The mean hitting time of a Markov chain is then defined for any $i\in\oneton$ by
\begin{align*}
\mathcal{M}(P)=\sum_{j=1}^n\bm{\pi}_{j}m_{ij}
\end{align*}
Due to its independence of the initial state $i$, the mean hitting time $\mathcal{M}(P)$ is also known as the Kemeny's constant, and it has been extensively studied in the literature \cite{JJH:13,SK:14}. We can also write the mean hitting time as
\begin{align*}
\mathcal{M}(P)=\sum_{i=1}^n\bm{\pi}_{i}\sum_{j=1}^n\bm{\pi}_{j}m_{ij}=\sum_{i=1}^n\sum_{j=1}^n\bm{\pi}_{i}\bm{\pi}_{j}m_{ij},
\end{align*}
which clearly shows that $\mathcal{M}(P)$ measures the expected time it takes for the surveillance agent to transition between two locations randomly selected according to the stationary distribution. The mean hitting time directly relates to the eigenvalues of the transition matrix as follows \cite[Eq. (9)]{JJH:13}: let $\lambda_1=1$ and $|\lambda_i|\leq1$ for $i\in\{2,\dots,n\}$ be the eigenvalues of an irreducible transition matrix $P$, then 
\begin{equation}\label{eq:KemenyEigen}
\mathcal{M}(P)=1+\sum_{j=2}^n\frac{1}{1-\lambda_j}.
\end{equation}

In robotic surveillance applications, if we require the surveillance robot to visit different locations as fast as possible, then we could adopt a surveillance strategy that minimizes the mean hitting time, and we obtain it by solving the following optimization problem.
\begin{problem}[Mean hitting time minimization for general Markov chains]\label{prob:meanhitting}
\emph{Given a strongly connected directed graph $\mathcal{G}=(V,\mathcal{E})$ and a stationary distribution $\bm{\pi}$, find a Markov chain with the stationary distribution $\bm{\pi}$ that has the minimum mean hitting time, i.e., solve the following optimization problem:}
\begin{equation*}
\begin{aligned}
& \underset{P\in\mathbb{R}^{n\times n}}{\textup{minimize}}
& & \mathcal{M}(P) \\
& \textup{subject to}
& & P\mathbb{1}_n = \mathbb{1}_n,\\
&&&\bm{\pi}^\top P=\bm{\pi}^\top,\\
&&& p_{ij}\geq0\quad\textup{if }(i,j)\in\mathcal{E},\\
&&&  p_{ij}=0\quad\textup{if }(i,j)\notin\mathcal{E}.
\end{aligned}
\end{equation*}
\end{problem}

\begin{figure}[http]
  \centering
  \subfigure[Pixel image of $P$]{\includegraphics[width=0.4\textwidth]{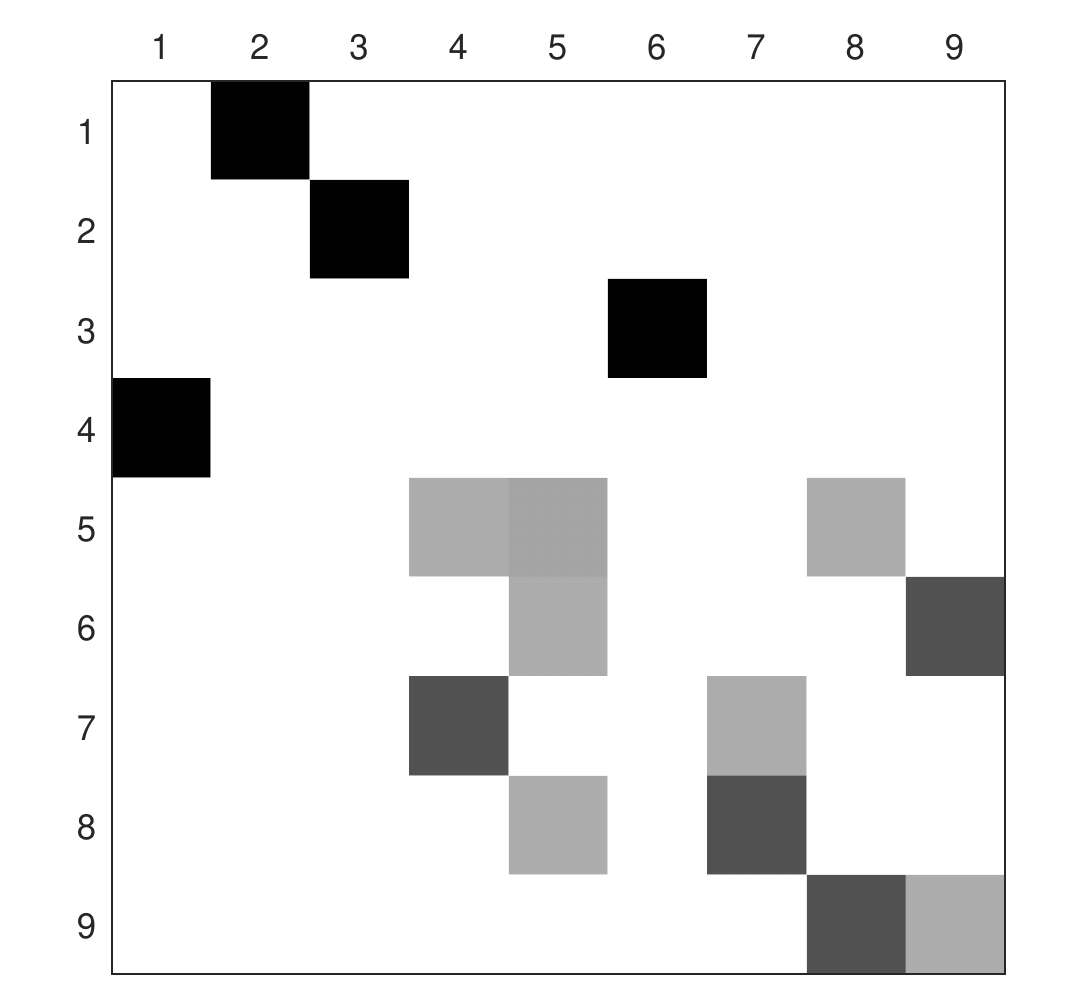}}\hfil
  \subfigure[Graphical image of $P$]{\includegraphics[width=0.35\textwidth]{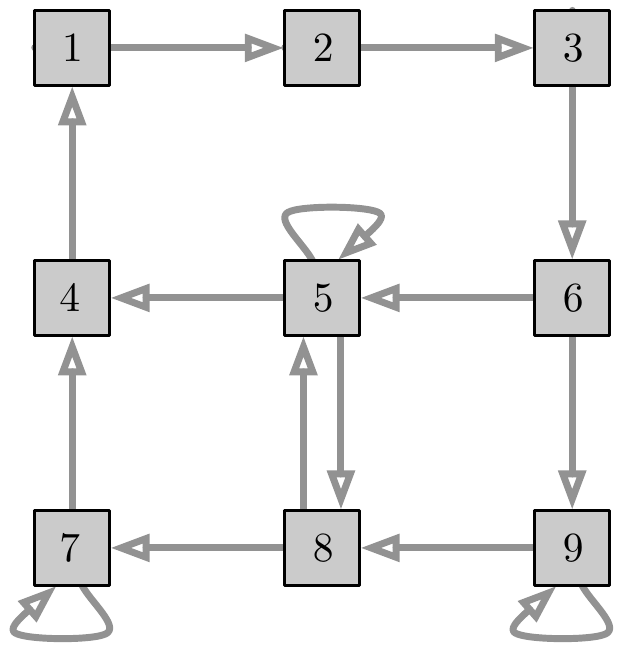}}
  \caption{Optimal strategy with minimum mean hitting time $\mathcal{M}(P)=6.78$ corresponding to Problem~\ref{prob:meanhitting} in a $3\times3$ grid.}\label{fig:minKemenygrid}
\end{figure}

Problem~\ref{prob:meanhitting} is a difficult nonconvex optimization problem, and the optimal solution is only available in special cases~\cite{SK:14}. For an illustration, we solve Problem~\ref{prob:meanhitting} for the grid graph using fmincon function in MATLAB. Since the problem is nonconvex, we start from $100$ different initial conditions and pick the best to be the optimal solution. A pixel image of the obtained solution is shown in Fig.~\ref{fig:minKemenygrid}, where each row has the transition probabilities from the corresponding location to other locations and the darkness of a pixel is proportional to the magnitude of the transition probability. One prominent feature of the optimal strategy is sparsity. Intuitively, it can be thought of approximately as a Hamiltonian tour that also satisfies the visit frequency.

Problem~\ref{prob:meanhitting} becomes tractable if we restrict the set of candidate Markov chains to be  reversible, in which case the mean hitting time is a convex function of the transition matrix in the domain. In fact, by using~\eqref{eq:KemenyEigen}, we can write the mean hitting time as \cite[Theorem 2]{RP-PA-FB:14b}
\begin{equation*}
\mathcal{M}(P)=\textup{Tr}((I_n-\diag(\bm{\pi})^{\frac{1}{2}}P\diag(\bm{\pi})^{-\frac{1}{2}}+\sqrt{\bm{\pi}}\sqrt{\bm{\pi}}^\top)^{-1}),
\end{equation*}
where $\textup{Tr}(\cdot)$ is the trace of a matrix and the square root is elementwise. Based on this, we formulate a convex optimization problem and solve for the Markov strategy with the minimum mean hitting time.
\begin{problem}[Mean hitting time minimization for reversible Markov chains]\label{prob:meanhittingrev}
\emph{Given a strongly connected directed graph $\mathcal{G}=(V,\mathcal{E})$ and a stationary distribution $\bm{\pi}$, find a reversible Markov chain with the stationary distribution $\bm{\pi}$ that has the minimum mean hitting time, i.e., solve the following optimization problem:}
	\begin{equation*}
	\begin{aligned}
	& \underset{P\in\mathbb{R}^{n\times n}}{\textup{minimize}}
	& & \textup{Tr}((I_n-\diag(\bm{\pi})^{\frac{1}{2}}P\diag(\bm{\pi})^{-\frac{1}{2}}+\sqrt{\bm{\pi}}\sqrt{\bm{\pi}}^\top)^{-1}) \\
	& \textup{subject to}
	& & P\mathbb{1}_n = \mathbb{1}_n,\\
	&&&\bm{\pi}_i p_{ij}=\bm{\pi}_jp_{ji}\quad\textup{for all } (i,j)\in \mathcal{E},\\
	&&& p_{ij}\geq0\quad\textup{if }(i,j)\in\mathcal{E},\\
	&&&  p_{ij}=0\quad\textup{if }(i,j)\notin\mathcal{E}.
	\end{aligned}
	\end{equation*}
\end{problem}

\begin{figure}[http]
  \centering
    \subfigure[Pixel image of $P$]{\includegraphics[width=0.4\textwidth]{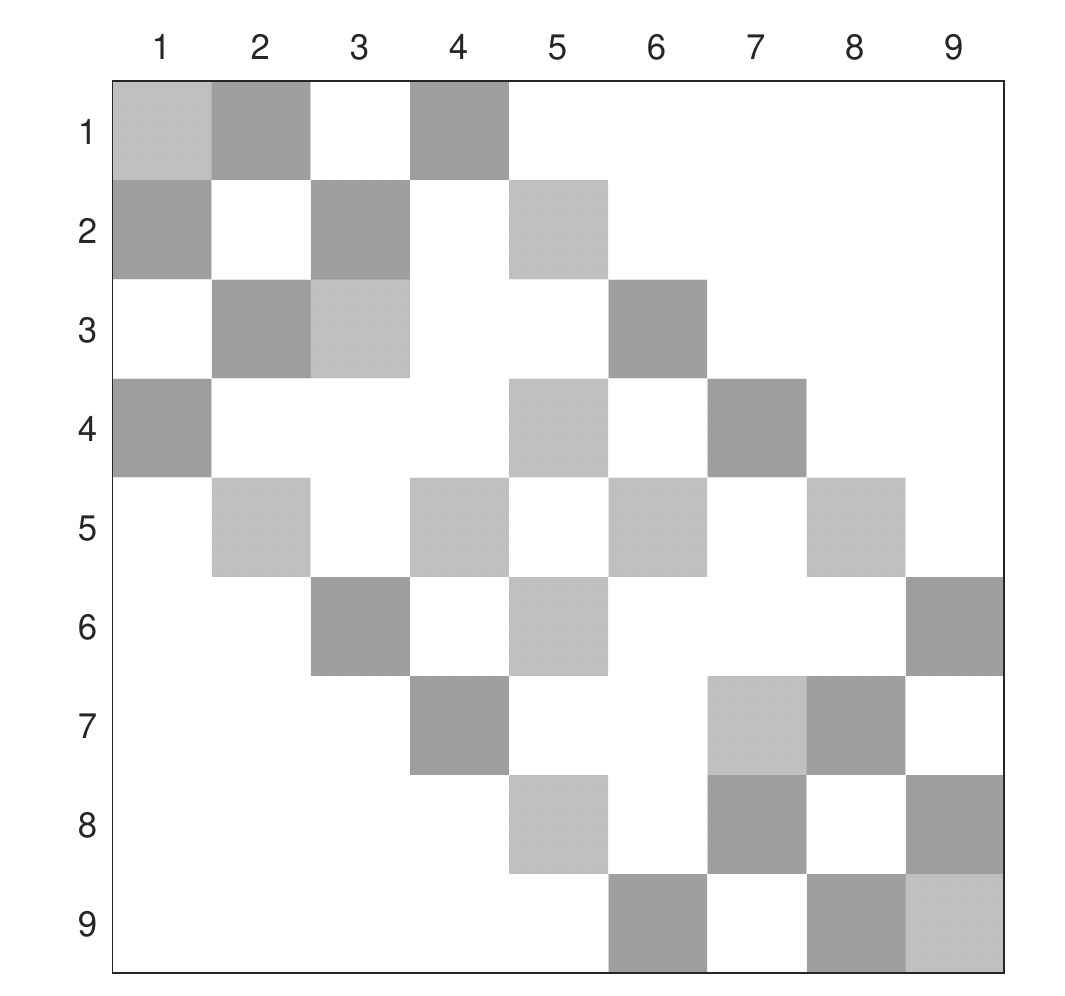}}\hfil
    \subfigure[Graphical image of $P$]{\includegraphics[width=0.35\textwidth]{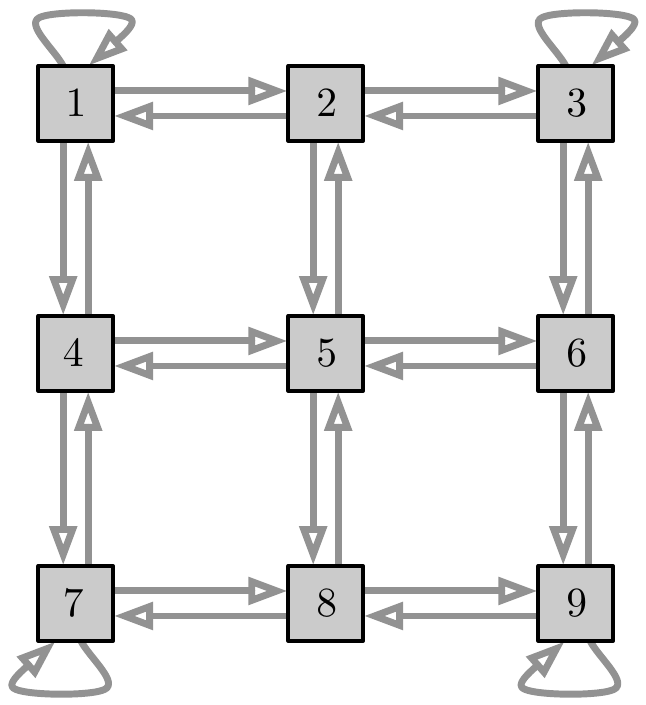}}
    \caption{Optimal reversible strategy with minimum mean hitting time $\mathcal{M}(P)=12.43$ corresponding to Problem~\ref{prob:meanhittingrev} in a $3\times3$ grid.}\label{fig:minKemenygridrev}
\end{figure}

The convexity of the objective function in Problem~\ref{prob:meanhittingrev} follows from the fact that the trace of the inverse of a positive definite matrix is convex~\cite{AG-SB-AS:08}. Moreover, Problem~\ref{prob:meanhittingrev} has a semidefinite reformulation and can be solved efficiently~\cite[Problem 2]{RP-PA-FB:14b}. Note that in both Problem~\ref{prob:meanhitting} and \ref{prob:meanhittingrev}, although we do not impose combinatorial irreducibility constraints on the Markov chain, the resulting solution must be irreducible as long as the constraint set is feasible and the given graph is strongly connected. This is because a reducible solution leads to infinite mean hitting time. The Markov chain with minimum mean hitting time is applied to the quickest detection of environmental anomalies in~\cite{PA-FB:15e}. As a comparison, we solve Problem~\ref{prob:meanhittingrev} using CVX~\cite{MG-SB:11-cvx} and plot the optimal solution in Fig.~\ref{fig:minKemenygridrev}. The optimal strategy given in Fig.~\ref{fig:minKemenygridrev} is denser than that given in Fig.~\ref{fig:minKemenygrid}. Moreover, the Kemeny's constant is almost doubled in the reversible case. 

\subsection{Arbitrary Travel Times}\label{subsec:hittingarb}
When there are travel times between neighboring locations in the graph, we can define the weighted mean hitting time similar to the unweighted case in Subsection~\ref{subsec:hittingunit}. Let $m_{ij}^{\textup{w}}$ be the expectation of the hitting time $T_{ij}^{\textup{w}}$ defined in~\eqref{eq:hittingtimeweight}, then
\begin{equation}\label{eq:meanhitttingtimeweight}
m_{ij}^{\textup{w}}=p_{ij}w_{ij} + \sum_{k\neq j}p_{ik}(w_{ik}+m_{kj}^{\textup{w}})=\sum_{k=1}^np_{ik}w_{ik}+\sum_{k=1}^np_{ik}m_{kj}^{\textup{w}}-p_{ij}m_{jj}^{\textup{w}}.
\end{equation}
We can write~\eqref{eq:meanhitttingtimeweight}  in matrix from as
\begin{equation}\label{eq:meanhittingmatrixweighted}
M^{\textup{w}}=(P\circ W)\mathbb{1}_n\mathbb{1}_n^\top+P(M^{\textup{w}}-\diag(M^{\textup{w}})),
\end{equation}
where $\circ$ is the elementwise product. By left multiplying the stationary distribution $\bm{\pi}^\top$ on both sides of~\eqref{eq:meanhittingmatrixweighted}, we have that
\begin{equation}\label{eq:expectedreturn}
\diag(M^{\textup{w}})=(\bm{\pi}^\top(P\circ W)\mathbb{1}_n)\cdot\diag(\bm{\pi})^{-1},
\end{equation}
which implies that the mean return time to location $i$ is still inversely proportional to $\bm{\pi}_i$ as in the unweighted case. Different from the case of unit travel times, the weighted sum $\sum_{j=1}^n\bm{\pi}_{j}m^{\text{w}}_{ij}$ depends on the initial state $i$ and the weighted mean hitting time is then defined as
\begin{align*}
\mathcal{M}^{\textup{w}}(P)=\sum_{i=1}^n\sum_{j=1}^n\bm{\pi}_{i}\bm{\pi}_{j}m^{\text{w}}_{ij}.
\end{align*}
It can be shown that the following relationship between the weighted mean hitting time and the mean hitting time
holds true~\cite[Theorem 8]{RP-PA-FB:14b}
\begin{align*}
\mathcal{M}^{\textup{w}}(P)=(\bm{\pi}^\top(P\circ W)\mathbb{1}_n)\cdot\mathcal{M}(P),
\end{align*}
where the impact of the travel times is factored out. To obtain a fast Markov chain over a weighted graph efficiently, we still restrict the set of candidate Markov chains to be reversible and formulate the following optimization problem.
\begin{problem}[Weighted mean hitting time minimization for Markov chains]\label{prob:meanhittingrevweighted}
\emph{Given a strongly connected weighted directed graph $\mathcal{G}=(V,\mathcal{E},W)$ and a stationary distribution $\bm{\pi}$, find a reversible Markov chain with the stationary distribution $\bm{\pi}$ that has the minimum weighted mean hitting time, i.e., solve the following optimization problem:}
	\begin{equation*}
	\begin{aligned}
	& \underset{P\in\mathbb{R}^{n\times n}}{\textup{minimize}}
	& & (\bm{\pi}^\top(P\circ W)\mathbb{1}_n)\cdot\textup{Tr}((I_n-\diag(\bm{\pi})^{\frac{1}{2}}P\diag(\bm{\pi})^{-\frac{1}{2}}+\sqrt{\bm{\pi}}\sqrt{\bm{\pi}}^\top)^{-1}) \\
	& \textup{subject to}
	& & P\mathbb{1}_n = \mathbb{1}_n,\\
	&&&\bm{\pi}_i p_{ij}=\bm{\pi}_jp_{ji}\quad\textup{for all } (i,j)\in \mathcal{E},\\
	&&& p_{ij}\geq0\quad\textup{if }(i,j)\in\mathcal{E},\\
	&&&  p_{ij}=0\quad\textup{if }(i,j)\notin\mathcal{E}.
	\end{aligned}
	\end{equation*}
\end{problem}

\begin{figure}[http]
	\centering{
	  \subfigure[Nonreversible strategy]{	
	\includegraphics[width=0.4\textwidth]{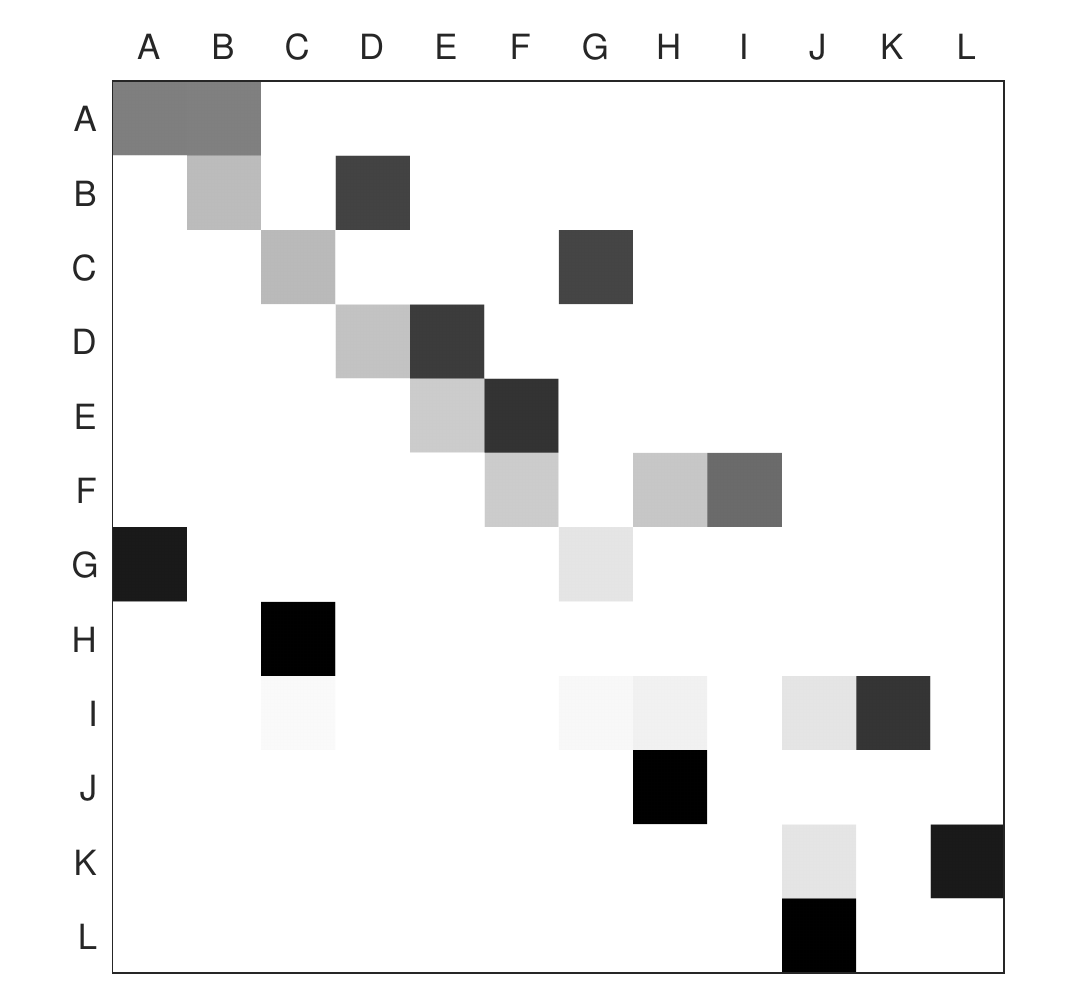}} \hfil
	\subfigure[Reversible strategy]{		
	\includegraphics[width=0.4\textwidth]{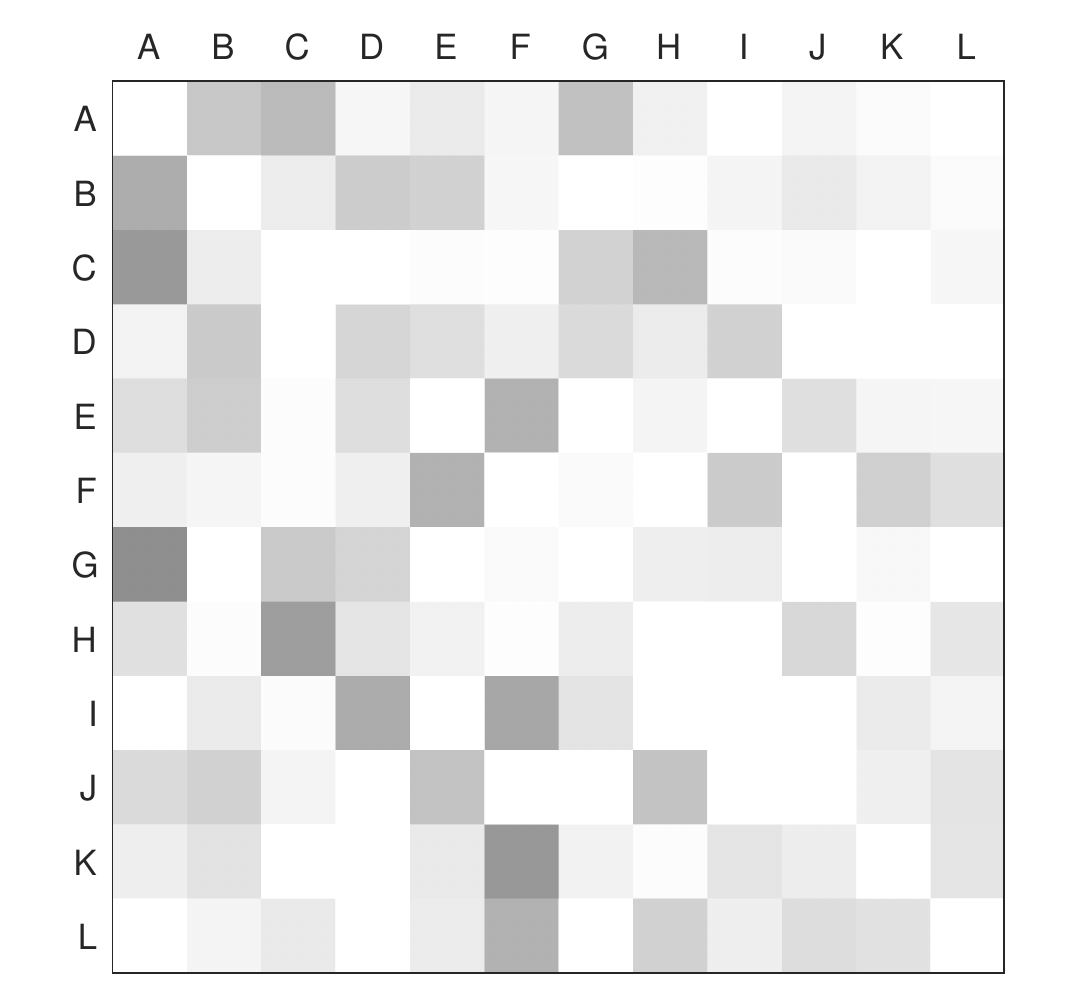}}}
\caption{Optimal Strategies with minimum mean hitting time in the SF map: (a) nonreversible solution with $\mathcal{M}^{\textup{w}}(P)=16.06$; (b) reversible solution $\mathcal{M}^{\textup{w}}(P)=29.88$ corresponding to Problem~\ref{prob:meanhittingrevweighted}.} \label{fig:SFKemeny}
\end{figure}

It turns out that by a change of variable, Problem~\ref{prob:meanhittingrevweighted} can also be recast to an equivalent semidefinite programming and thus solved efficiently. We solve Problem~\ref{prob:meanhittingrevweighted} and its nonreversible version for the SF map data set. The solutions are illustrated in Fig.~\ref{fig:SFKemeny}. Similar to the unweighted cases discussed in Subsection~\ref{subsec:hittingunit}, the optimal solution to Problem~\ref{prob:meanhittingrevweighted} is denser than its nonreversible counterpart and has a higher (worse) optimal value. 

It is worth mentioning that the computational tractability of Problems~\ref{prob:meanhittingrev} and~\ref{prob:meanhittingrevweighted} relies critically on the reversibility condition on the set of candidate Markov chains. However, reversible Markov chains are known to be slow~\cite{KJ-DS-JS:10}, which is also illustrated by our numerical examples. Therefore, how to solve Problem~\ref{prob:meanhitting} and its weighted version efficiently with optimality guarantees remains an open and interesting problem.

\subsection{Hitting Times for Multiple Surveillance Robots}\label{subsec:multihitting}
In this subsection, we consider a team of robots indexed by $\{1,\dots,N\}$, and each robot $i$ patrols a strongly connected graph $\mathcal{G}^i=(V,\mathcal{E}^{i})$ according to an irreducible Markov strategy $P^i$. Note that the patrolling graphs have the common set of nodes but possibly different sets of edges. The hitting time for this team of robots from their current locations $i_1,\dots i_N$ to a location $j$ is defined as the first time at least one of the robots arrives at location $j$, i.e.,
\begin{align}\label{eq:hittingtimemultiple}
\begin{split}
T_{i_1\dots i_N,j}=\min\{k\geq1\,|\,&X^{1}_k=j\textup{ or }X^{2}_k=j\cdots\textup{ or }X^{N}_k=j,\\
&X^{h}_0=i_h\textup{ for }h\in\{1,\dots,N\}\},
\end{split}
\end{align}
where $X_k^{h}$ is the location of robot $h$ at time $k$. Let $m_{i_1\dotsi_N,j}$ be the expectation of $T_{i_1\dotsi_N,j}$, then similar to the case of single robot before, we have that
\begin{align}\label{eq:multirobothitting}
\begin{split}
m_{i_1\dots i_N,j}&=1-\prod\limits_{h=1}^N (1-p^h_{i_hj}) + \sum_{k_1\neq j}\cdots\sum_{k_N\neq j}p^1_{i_1k_1}\cdots p^N_{i_Nk_N}(1+m_{k_1\dots k_N,j})\\
&=1+\sum_{k_1\neq j}\cdots\sum_{k_N\neq j}p^1_{i_1k_1}\cdots p^N_{i_Nk_N}m_{k_1\dots k_N,j}.
\end{split}
\end{align}
In order to write~\eqref{eq:multirobothitting} in matrix form, we collect $m_{i_1\dots i_N,j}$ into $M\in\mathbb{R}^{n^N\times n}$ where each row of $M$ contains mean hitting times from a specific robot team ``configuration'' (a set of locations currently occupied by the team of robots) to different locations. The rows of $M$ are arranged by cycling through the possible locations of robots. Concretely, the mean hitting time matrix $M$ has the following structure
\begin{equation*}
M=\begin{bmatrix}
m_{1\cdots11,1}&m_{1\cdots11,2}&\cdots&m_{1\cdots11,n}\\
m_{1\cdots12,1}&m_{1\cdots12,2}&\cdots&m_{1\cdots12,n}\\
\vdots&\vdots&\vdots&\vdots\\
m_{1\cdots1 n,1}&m_{1\cdots1 n,2}&\cdots&m_{1\cdots1 n,n}\\
m_{1\cdots21,1}&m_{1\cdots21,2}&\cdots&m_{1\cdots21,n}\\
\vdots&\vdots&\vdots&\vdots\\
m_{n\cdots nn,1}&m_{n\cdots nn,2}&\cdots&m_{n\cdots nn,n}\\
\end{bmatrix}.
\end{equation*}
Given the structure of $M$, we can write \eqref{eq:multirobothitting} in matrix form as \cite[Lemma 3.1]{RP-AC-FB:14k}
\begin{equation}\label{eq:multirobothittingmatrix}
M=\mathbb{1}_{n^N}\mathbb{1}_{n}^\top+(P^1\otimes\cdots\otimes P^N)(M-M_d),
\end{equation}
where $\otimes$ is the Kronecker product and $M_d\in\mathbb{R}^{n^N\times n}$ has nonzero elements $m_{i_1\dots i_N,j}$ only at locations $(i_1\dots i_N,j)$ that satisfy $i_h=j$ for at least one $h\in\{1,\dots,N\}$. To calculate the mean hitting times, we vectorize both sides of~\eqref{eq:multirobothittingmatrix} and obtain
\begin{equation}\label{eq:multiplemeanhitting}
\vecz(M)=(I_{n^{N+1}}-(I_n\otimes P^1\otimes\cdots\otimes P^N)(I_{n^{N+1}}-E))^{-1}\mathbb{1}_{n^{N+1}},
\end{equation}
where $E\in\mathbb{R}^{n^{N+1}\times n^{N+1}}$ is a binary matrix that satisfies $E\vecz(M)=\vecz({M_d})$. The matrix inversion in~\eqref{eq:multiplemeanhitting} is well defined since the robots' surveillance strategies are all irreducible.

In the multi-robot case, the problem dimension grows exponentially fast when we scale up the number of surveillance robots. Therefore, it is not obvious how the strategies for a team of robots could be jointly optimized in an efficient manner. One can certainly define a mean hitting time for the team of robots analogous to the single robot case. However, the interpretation is not as clean, and the optimization simply becomes intractable. Designing multi-robot surveillance strategies as a whole is still an ongoing problem.

\subsection{Meeting Times of Two Markovian Agents}\label{subsec:meetingtime}
In this subsection, we study the meeting times of two Markovian agents. The development still has the hitting times as the core elements and also has similarity to the case of multiple surveillance agents. We consider a surveillance robot with a Markov strategy $P^\textup{p}$ and an evader with a Markov strategy $P^\textup{e}$, both defined over a common strongly connected directed graph $\mathcal{G}=(V,\mathcal{E})$. Starting from two initial positions $i,j\in V$, the meeting time of the two agents is defined by
\begin{equation*}
T_{ij}=\min\{k\geq1\,|\,X^{\textup{p}}_k=X_k^{\textup{e}},X^{\textup{p}}_0=i,X^{\textup{e}}_0=j\},
\end{equation*}
where $X_k^{\textup{p}}$  and $X_k^{\textup{e}}$ are the locations of the surveillance robot and the evader at time $k\geq0$, respectively. It should be emphasized that here, the subscripts $i$ and $j$ of $T_{ij}$ represent the initial locations of two different agents. If the evader is stationary whose strategy is given by $P^{\textup{e}}=I_n$, then the meeting time $T_{ij}$ corresponds to exactly the hitting time of the surveillance agent with the strategy $P^{\textup{p}}$. Let $m_{ij}$ be the expectation of $T_{ij}$, then it follows that
\begin{equation}\label{eq:meetingtime}
m_{ij}  =\sum_{k=1}^n p_{ik}^{\textup{p}}p_{jk}^{\textup{e}} + \sum_{k_1 \neq h_1 } p_{ik_1}^{\textup{p}}p_{jh_1}^{\textup{e}} (1+m_{k_1h_1})=1+\sum_{k_1 \neq h_1 } p_{ik_1}^{\textup{p}}p_{jh_1}^{\textup{e}} m_{k_1h_1}.
\end{equation}
Once again, we can organize $m_{ij}$ into a matrix $M\in\mathbb{R}^{n\times n}$ and obtain
\begin{equation}\label{eq:meetingtimematrix}
M = \mathbb{1}_n \mathbb{1}_n^\top + P^{\textup{p}}(M-\diag(M)){P^{\textup{e}}}^\top.
\end{equation}

\begin{figure}[http]
	\centering	
	\subfigure[Finite meeting time]{			
		\begin{tikzpicture}
		\node[] at (-1.3,0) {\large$P^{\textup{p}}$};
		\node[state,minimum size = 0.3cm] at (0, 0) (nodeone) {\text{$1$}};
		\node[state,minimum size = 0.3cm] at (2,0)     (nodetwo)     {\text{$2$}};
		\node[] at (-1.3,-1) {\large$P^{\textup{e}}$};
		\node[state,minimum size = 0.3cm] at (0, -1) (nodethree) {\text{$1$}};
		\node[state,minimum size = 0.3cm] at (2,-1)     (nodefour)     {\text{$2$}};
		\draw[every loop,
		auto=right,
		>=latex,
		]
		(nodetwo)     edge[bend left=15, auto=right] node {} (nodeone)	
		(nodeone) edge[bend left=15, auto=right] node {} (nodetwo)	
		(nodefour)     edge[bend left=15, auto=right] node {} (nodethree)	
		(nodethree) edge[bend left=15, auto=right] node {} (nodefour)
		(nodethree) edge[loop left,auto=right] node {} (nodethree)	
		(nodefour) edge[loop right,auto=right] node {} (nodefour)	
		(nodeone) edge[loop left,auto=right] node {} (nodeone)	
		(nodetwo) edge[loop right,auto=right] node {} (nodetwo)
		
		;
		\end{tikzpicture}
	}	
	\hfil
	\subfigure[Infinite meeting time]{			
	\begin{tikzpicture}
	\node[] at (-1,0) {\large$P^{\textup{p}}$};
	\node[state,minimum size = 0.3cm] at (0, 0) (nodeone) {\text{$1$}};
	\node[state,minimum size = 0.3cm] at (2,0)     (nodetwo)     {\text{$2$}};
	\node[] at (-1,-1) {\large$P^{\textup{e}}$};
	\node[state,minimum size = 0.3cm] at (0, -1) (nodethree) {\text{$1$}};
	\node[state,minimum size = 0.3cm] at (2,-1)     (nodefour)     {\text{$2$}};
	\draw[every loop,
	auto=right,
	>=latex,
	]
	(nodetwo)     edge[bend left=15, auto=right] node {} (nodeone)	
	(nodeone) edge[bend left=15, auto=right] node {} (nodetwo)	
	(nodefour)     edge[bend left=15, auto=right] node {} (nodethree)	
	(nodethree) edge[bend left=15, auto=right] node {} (nodefour);
	\end{tikzpicture}
}
	\caption{The pursuer and evader in (a) have finite meeting times as starting from any positions for the agents, there exist paths to the common locations $(1,1)$ and $(2,2)$. However, in (b) there exist no paths to common locations from locations $(1,2)$ and $(2,1)$.}\label{fig:twonodeexample}
\end{figure}
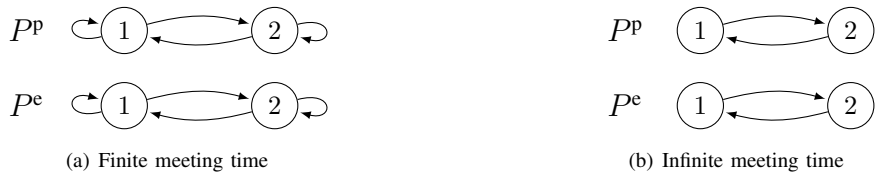

Unlike previous cases where irreducibility of individual strategies is sufficient to ensure finiteness of mean hitting times, the mean meeting time may be infinite even if both $P^{\textup{p}}$ and $P^{\textup{e}}$ are irreducible. In fact, in order for two Markovian agents to meet in finite time from any initial positions, there must exist paths of equal length that lead a common position for the agents from their respective initial positions; see a simple example in Fig.~\ref{fig:twonodeexample} for an illustration. The finiteness of meeting times can be studied systematically by looking at the Kronecker product of graphs~\cite{PMW:62}, and we refer the interested readers to~\cite{XD-MG-RP-FB:14l} for more details.

With the mean meeting time matrix $M$ and the stationary distributions $\bm{\pi}^{\textup{p}}$ and $\bm{\pi}^{\textup{e}}$ for the surveillance agent and the evader, we can define the expected meeting time as
\begin{equation}\label{eq:meanmeetingtime}
\mathcal{M}={\bm{\pi}^{\textup{p}}}^\top M \bm{\pi}^{\textup{e}}.
\end{equation}
Since we do not require the strategies for the surveillance agent and the evader to be irreducible, the stationary distributions $\bm{\pi}^{\textup{p}}$ and $\bm{\pi}^{\textup{e}}$ may not be unique. However, the expected meeting time in~\eqref{eq:meanmeetingtime} can still be defined as long as all entries of $M$ are finite and one set of stationary distributions are specified. Given the evader's strategy $P^\textup{e}$ and a stationary distribution $\bm{\pi}^{\textup{e}}$, we can solve the following optimization problem to find a surveillance strategy $P^{\textup{p}}$ that captures (meets) the evader as fast as possible.
\begin{problem}[Expected meeting time minimization]\label{prob:meanmeetingtime}
\emph{Given a strongly connected directed graph $\mathcal{G}=(V,\mathcal{E})$, an evader's strategy $P^\textup{e}$ with a stationary distribution $\bm{\pi}^\textup{e}$, and a stationary distribution $\bm{\pi}^\textup{p}$ for the surveillance agent. Find a Markov chain with the stationary distribution $\bm{\pi}^\textup{p}$ that has the minimum expected meeting time, i.e., solve the following optimization problem:}
 	\begin{equation*}
	\begin{aligned}
	&\underset{{P^\textup{p}}\in\mathbb{R}^{n\times n}}{\textup{minimize}}
	& &{\bm{\pi}^{\textup{p}}}^\top M \bm{\pi}^{\textup{e}}\\
	& \textup{subject to}
	&& {\bm{\pi}^\textup{p}}^\top P^\textup{p}={\bm{\pi}^\textup{p}}^\top,\\
	&&&P^\textup{p}\mathbb{1}_n=\mathbb{1}_n,\\
	&&&p_{i,j}^{\textup{p}}\geq0\quad \textup{if }(i,j)\in \mathcal{E},\\
	&&& p_{i,j}^{\textup{p}}=0\quad \textup{if }(i,j)\notin \mathcal{E}.
	\end{aligned}
	\end{equation*}	
\end{problem}

\begin{figure}[http]
  \centering
  \subfigure[Pixel image of $P$]{\includegraphics[width=0.4\textwidth]{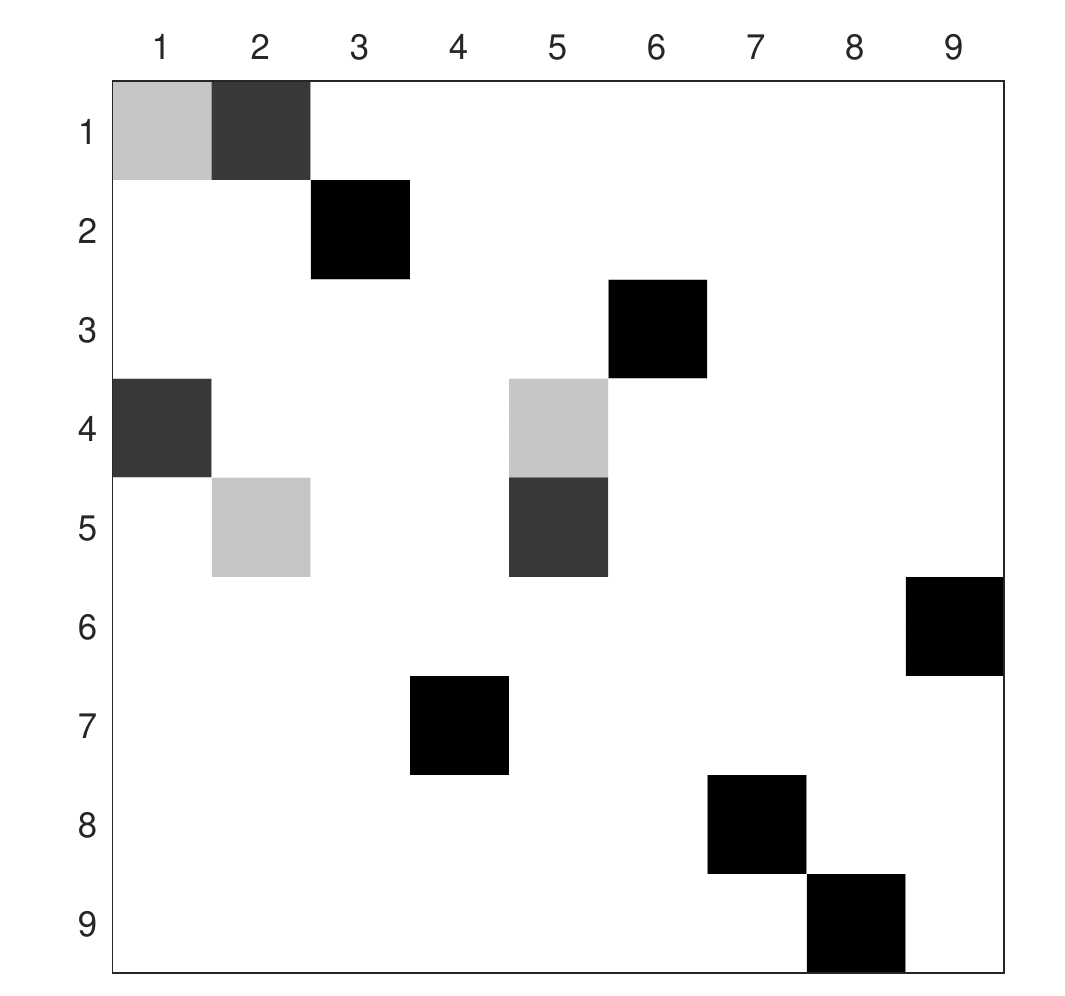}}\hfil
  \subfigure[Graphical image of $P$]{\includegraphics[width=0.35\textwidth]{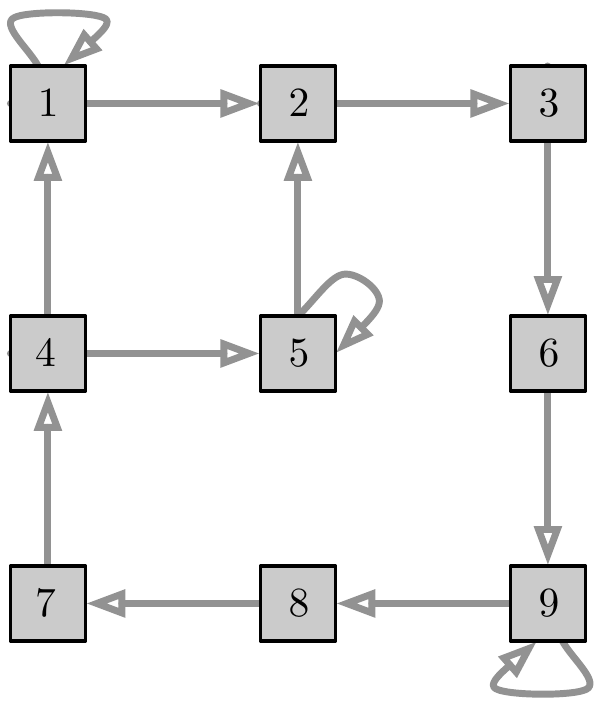}}
  \caption{Optimal strategy with minimum expected meeting time $\mathcal{M}(P)=9.86$ against a randomly walking evader corresponding to Problem~\ref{prob:meanmeetingtime} in a $3\times3$ grid.}\label{fig:minmeeting}
\end{figure}

The expected meeting time minimization Problem~\ref{prob:meanmeetingtime} has the mean hitting time minimization Problem~\ref{prob:meanhitting} as a special case and thus is also difficult to solve. In fact, we can interpret the mean hitting time minimization problem as an expected meeting time minimization problem: a Markovian surveillance agent is chasing a stationary evader whose strategy is given by $P^{\textup{e}}=I_n$ and $\bm{\pi}^\textup{e}=\bm{\pi}^\textup{p}$. Other than the special cases discussed in~\cite{XD-MG-RP-FB:14l}, we have to resort to numerical solvers in order to obtain a locally optimal solution. We plot the optimal strategy for a pursuer chasing a randomly walking evader in Fig.~\ref{fig:minmeeting}, where the evader moves from its current location to neighboring locations (including a self loop) with equal probabilities. The pursuer's strategy has a strong sparsity pattern similar to that shown in Fig.~\ref{fig:minKemenygrid}. Note that in this case, the pursuer has a uniform stationary distribution while the evader's stationary distribution is proportional to the node degrees.

\section{UNPREDICTABLE SURVEILLANCE: MAXIMIZATION OF ENTROPY}\label{sec:unpredictableMC}
One of the motivations for adopting Markov chains as surveillance strategies is to introduce stochasticity to the motion plans of the surveillance robot. In this section, we design surveillance strategies that are maximally unpredictable and thus hard for potential intruders to exploit. We utilize two different notions of entropy to quantify unpredictability of Markov chains. The optimization of the entropy rate in Subsection~\ref{subsec:entropyrate} is studied in~\cite{MG-SJ-FB:17b} and the treatment of the return time entropy in Subsection~\ref{subsec:returnentropy} follows from \cite{XD-MG-FB:17o}.

\subsection{Maximization of Entropy Rate}\label{subsec:entropyrate}
As introduced in Section~\ref{sec:entropyrateintro}, the entropy rate is a classic metric that measures the randomness of the random sequence generated by Markov chains. In robotic surveillance, this sequence is composed of locations visited by the surveillance robot. Since the travel times on the graph do not affect the entropy rate of the strategy, we consider unweighted graphs in this subsection. Moreover, we allow self loops at all locations and assume the graph is symmetric as it turns out that the optimal strategy in this setting can be computed very efficiently. In order to find a Markov chain that maximizes the entropy rate defined in~\eqref{eq:entropyrate}, we formulate the following optimization problem.
\begin{problem}[Entropy rate maximization]\label{prob:entropyrate}
\emph{Given a strongly connected directed graph $\mathcal{G}=(V,\mathcal{E})$ with self loops and symmetric adjacency matrix and a stationary distribution $\bm{\pi}$, find a Markov chain with the stationary distribution $\bm{\pi}$ that has the maximum entropy rate, i.e., solve the following optimization problem:}
	\begin{equation*}
	\begin{aligned}
	& \underset{P\in\mathbb{R}^{n\times n}}{\textup{maximize}}
	& & \Hrate(P) \\
	& \textup{subject to}
	& & P\mathbb{1}_n = \mathbb{1}_n,\\
	&&&\bm{\pi}^\top P=\bm{\pi}^\top,\\
	&&& p_{ij}\geq0\quad\textup{if }(i,j)\in\mathcal{E},\\
	&&&  p_{ij}=0\quad\textup{if }(i,j)\notin\mathcal{E}.
	\end{aligned}
	\end{equation*}
\end{problem}

As in previous problems, a stationary distribution is specified in Problem~\ref{prob:entropyrate}. A similar problem where only the graph structure is given was considered in~\cite{ZB-JD-JML-BW:09}, where the optimal solution is carefully constructed. One may recognize that Problem~\ref{prob:entropyrate} is a convex problem even without requiring the graph to be symmetric and have self loops, and it can be solved by generic convex programming solvers. However, Problem~\ref{prob:entropyrate} exhibits several intriguing properties and allows for faster iterative algorithms. Let $A\in\mathbb{R}^{n\times n}$ be the unit-diagonal symmetric adjacency matrix of the strongly connected graph $\mathcal{G}=(V,\mathcal{E})$. We define the \emph{maxentropic vector map} $\phi:\mathbb{R}_{>0}^n\mapsto\mathbb{R}^n_{>0}$ by
\begin{equation*}
\phi(x)=\diag(x)Ax.
\end{equation*}
Moreover, we define the \emph{maxentropic matrix map} $\Phi:\mathbb{R}_{>0}^n\mapsto\mathbb{R}^{n\times n}_{\geq0}$ by
\begin{equation*}
\Phi(x)=\diag(Ax)^{-1}A\diag(x).
\end{equation*}
These two maps are essential for constructing the optimal solution to Problem~\ref{prob:entropyrate}, and many nice properties of them are summarized in~\cite[Theorem 3]{MG-SJ-FB:17b} and \cite[Theorem 2]{MG-SJ-FB:17b}. To solve Problem~\ref{prob:entropyrate}, we take the following steps:
\begin{enumerate}
\item\label{itm:invermap} find $x^*\in\mathbb{R}^n_{>0}$ such that $\phi(x^*)=\bm{\pi}$, i.e., solve the inverse maxentropic vector map $x^*=\phi^{-1}(\bm{\pi})$;
\item\label{itm:optimalsolution} construct the optimal solution to Problem~\ref{prob:entropyrate} as $P^*=\Phi(x^*)=\diag(Ax^*)^{-1}A\diag(x^*)$;
\item\label{itm:optimalvalue} compute the optimal value by $\Hrate(P^*)=-2{x^*}^\top A\diag(x^*)\log x^*+\bm{\pi}^\top\log\bm{\pi}$, where the logarithm is taken elementwise. 
\end{enumerate}
In the procedures outlined above, once the inverse maxentropic vector map $x^*=\phi^{-1}(\bm{\pi})$ in step~\ref{itm:invermap}  is solved, the optimal solution and its value can be obtained straightforwardly in steps~\ref{itm:optimalsolution} and~\ref{itm:optimalvalue}. 
First of all, as ensured by one of the properties of the inverse maxentropic vector map, step~\ref{itm:invermap} is well defined in the sense that for any $\bm{\pi}\in\mathbb{R}_{>0}^n$, there exists a unique $x^*$ that satisfies $\phi(x^*)=\bm{\pi}$. Moreover, let
\begin{equation*}
\eta=\max_{i}\{\sum_{j=1}^na_{ij}\sqrt{\bm{\pi}_j}\},\qquad x^0=\frac{1}{\sqrt{\max_{i}\{\sum_{j=1}^na_{ij}\bm{\pi}_j\}}}\bm{\pi},
\end{equation*}
then the following linear iteration converges to $x^*$ \cite[Theorem 4]{MG-SJ-FB:17b}
\begin{equation*}
x^{k+1}=x^k-\frac{1}{2\eta}(\diag(x^k)Ax^k-\bm{\pi}),~\textup{for }k\geq0.
\end{equation*}

\begin{figure}[http]
  \centering
  \subfigure[Pixel image of $P$]{\includegraphics[width=0.4\textwidth]{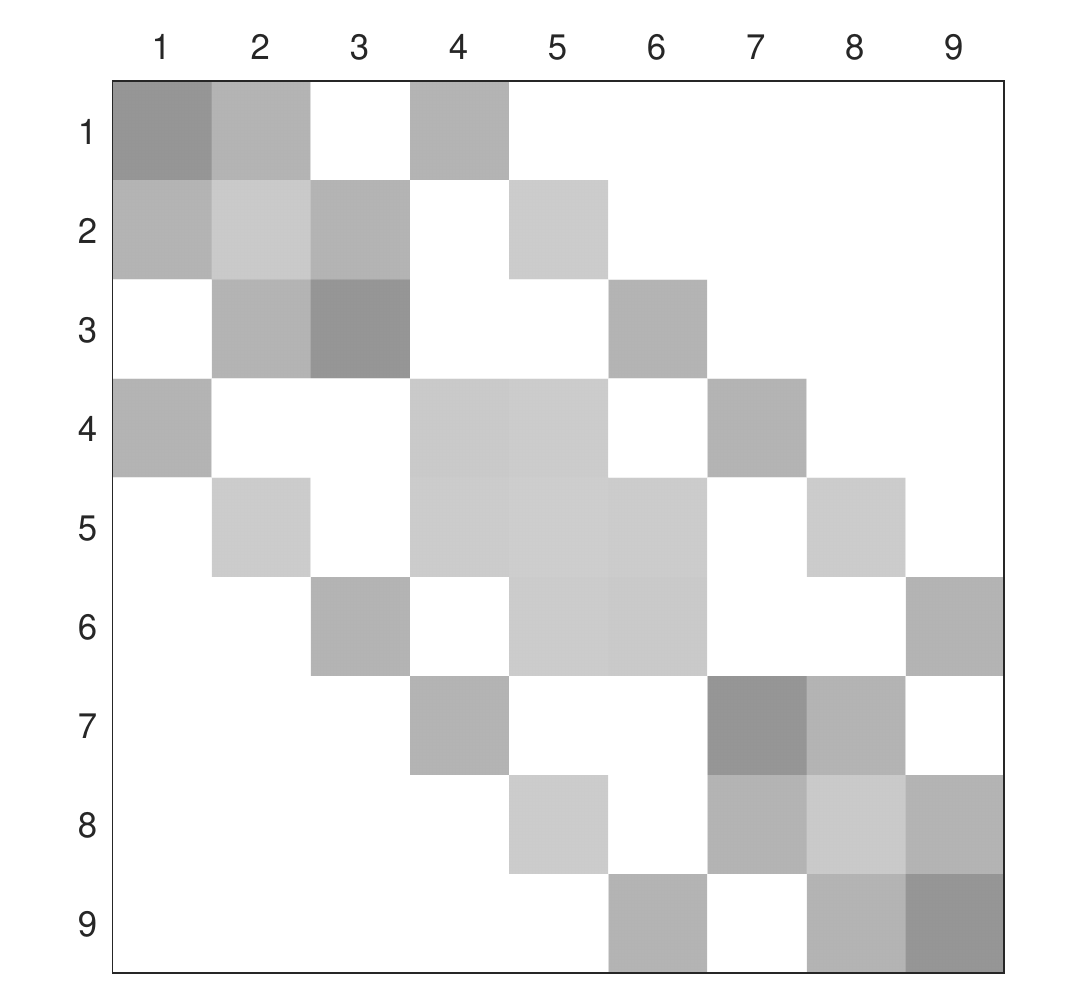}}\hfil
  \subfigure[Graphical image of $P$]{\includegraphics[width=0.35\textwidth]{graph-grid-topology}}  
  \caption{Optimal strategy with maximum entropy rate $\Hrate(P)=1.27$ corresponding to Problem~\ref{prob:entropyrate} in a $3\times3$ grid.}\label{fig:maxentropygrid}
\end{figure}

Numerical comparisons and computational complexity analysis are given in~\cite{MG-SJ-FB:17b} where the efficiency of the proposed linear iteration is clearly shown. It is also worth mentioning that the optimal solution to Problem~\ref{prob:entropyrate} is automatically reversible as proven in~\cite[Theorem 6]{MG-SJ-FB:17b}. Here, we solve Problem~\ref{prob:entropyrate} using the linear iteration and plot the solution in Fig.~\ref{fig:maxentropygrid}. The solution has a similar pattern as the reversible chain with minimum Kemeny's constant. However, the transition probabilities are shaped differently, and the strategy with maximum entropy rate has self loops at all locations. 

\subsection{Maximization of Return Time Entropy}\label{subsec:returnentropy}
The entropy rate is a ``global'' measure of unpredictability of Markov chains, i.e., it is concerned with the entire sequence of locations visited by the surveillance robot. However, from a local observer's point of view, the inter-visit times at a location of interest may be the only available information. For Markov chains, these inter-visit times correspond exactly to the first return times $T_{ii}$ for locations $i\in\oneton$. In order to make it hard for the observers (potential intruders) to predict when the surveillance robot will come back after it leaves the current location, the maximization of unpredictability in return times measured by the entropy is proposed. Note that the return times depend not only on the Markov strategy, but also on the travel times between locations. In this subsection, we consider strongly connected directed graphs with integer-valued edge weights. Assuming integer-valued travel times is not restrictive because we can always quantize and approximate real-valued travel times reasonably well by choosing appropriate quantization units. Moreover, it is much more convenient to characterize the return time distributions if the travel times are integer-valued. 

To ease the notation, we introduce the shorthand $F_k\in\mathbb{R}^{n\times n}$ where the $(i,j)$-th element of $F_k$ represents the probability that the hitting time $T_{ij}^\textup{w}$ takes value $k$, i.e., $F_{k}(i,j)=\mathbb{P}(T_{ij}^{\textup{w}}=k)$.  Then, we have the following recursion
\begin{equation}\label{eq:hittingprob}
F_k(i,j)=p_{ij}\mathbf{1}_{\{k=w_{ij}\}}+\sum_{h\neq j}p_{ih}F_{k-w_{ih}}(h,j),
\end{equation}
with the initial conditions $F_k=\mathbb{0}_{n\times n}$ for $k\leq0$. The formula~\eqref{eq:hittingprob} can also be arranged in a vectorized form as
\begin{equation}\label{eq:hittingtimeprob}
\vecz(F_k)=\vecz(P\circ \mathbf{1}_{\{k\mathbb{1}_n\mathbb{1}_n^\top=W\}}) +\sum_{i=1}^n\sum_{j=1}^np_{ij}(E_j\otimes\mathbb{e}_i\mathbb{e}_j^\top)\vecz(F_{k-w_{ij}}),
\end{equation}
where $\mathbf{1}_{\{k\mathbb{1}_n\mathbb{1}_n^\top=W\}}$ is a binary matrix whose $(i,j)$-th element is $1$ if $w_{ij}=k$ and $E_{i}=\diag(\mathbb{1}_n-\mathbb{e}_i)$. Despite its seemingly complicated form, the hitting time probabilities in~\eqref{eq:hittingtimeprob} actually satisfy a linear system with time delays and transition probabilities as inputs. Moreover, the inputs to this system are directly related to the system matrix, and the hitting time probabilities are polynomial functions of the transition probabilities.

Using the notation introduced above, we can rewrite the return time entropy of a Markov chain with travel times explicitly as
\begin{equation}\label{eq:returntimeentropyweight}
\Hrt(P)=-\sum_{i=1}^{n}\bm{\pi}_i\sum_{k=1}^\infty F_k(i,i)\log F_{k}(i,i).
\end{equation}
Our goal is to find a Markov chain that has the maximum return time entropy by solving the following optimization problem.
\begin{problem}[Return time entropy maximization]\label{prob:returntimeentropy}
\emph{Given a strongly connected directed weighted graph $\mathcal{G}=(V,\mathcal{E},W)$ and a stationary distribution $\bm{\pi}$, pick a minimum transition probability $\epsilon>0$, find a Markov chain with the stationary distribution $\bm{\pi}$ that has the maximum return time entropy, i.e., solve the following optimization problem:}
	\begin{equation*}
	\begin{aligned}
	& \underset{P\in\mathbb{R}^{n\times n}}{\textup{maximize}}
	& & \Hrt(P) \\
	& \textup{subject to}
	& & P\mathbb{1}_n = \mathbb{1}_n,\\
	&&&\bm{\pi}^\top P=\bm{\pi}^\top,\\
	&&& p_{ij}\geq\epsilon\quad\textup{if }(i,j)\in\mathcal{E},\\
	&&&  p_{ij}=0\quad\textup{if }(i,j)\notin\mathcal{E}.
	\end{aligned}
	\end{equation*}
\end{problem}

We clarify several technical details in Problem~\ref{prob:returntimeentropy} in the following:
\begin{enumerate}
\item the minimum transition probability $\epsilon$ ensures the irreducibility of Markov chains; moreover, this also ensures that the constraint set of Problem~\ref{prob:returntimeentropy} is a compact set of irreducible Markov chains.

\item the objective function $\Hrt(P)$ as given by~\eqref{eq:returntimeentropyweight} is an infinite series; fortunately, by analyzing the linear dynamical system~\eqref{eq:hittingtimeprob}, we can show that the return time distributions $F_k(i,i)$'s decay exponentially fast if $P$ is irreducible. Therefore, the infinite series is summable and well defined. In fact, the objective function continuously depends on the transition probabilities.
\end{enumerate}
The above clarifications imply that in Problem~\ref{prob:returntimeentropy}, a continuous function is optimized over a compact set. Therefore, there exists at least one optimal solution and Problem~\ref{prob:returntimeentropy} is well posed. Moreover, it should be pointed out that due to the stationary distribution constraint and by~\eqref{eq:expectedreturn}, the return times have bounded expected values and the surveillance robot will come back to each location within some time on average.

The difficulty of Problem~\ref{prob:returntimeentropy} partly comes from the fact that a closed form of $\Hrt(P)$ as a function of the transition probabilities is not available. However, in the particular case of a complete graph with unit travel times, the return times can be made to follow geometric distributions and their entropy has a closed form. Moreover, by the maximum entropy principle for discrete random variables~\cite{SG-AS:85}, the optimal solution is given by $P=\mathbb{1}_n\bm{\pi}^\top$.

As in~\cite{LE-TMC:93}, it is also possible to establish a connection between the return time entropy (in unweighted graphs) and the entropy rate as follows
\begin{equation}\label{eq:relation}
\Hrate(P)\leq\Hrt(P)\leq n\Hrate(P),
\end{equation}
where the inequalities are tight in the sense that they become equalities in certain cases. 
This relationship~\eqref{eq:relation} implies that $\Hrt(P)$ and $\Hrate(P)$ are two distinct metrics and are not interchangeable. In fact, they may differ by a factor of $n$, which is the dimension of the graph.

\begin{figure}[http]
	\centering
	\includegraphics[width=0.4\textwidth]{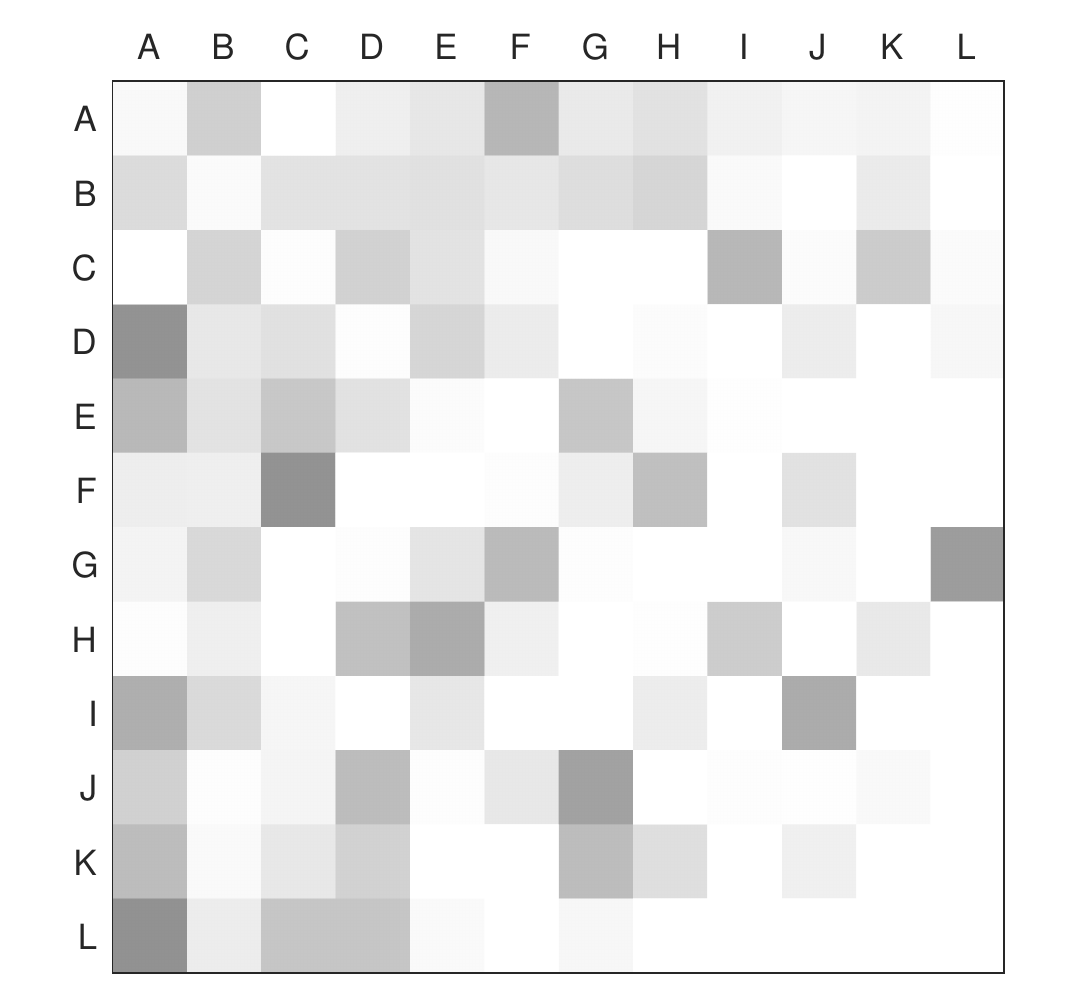}
	\caption{Optimal strategy with maximum return time entropy $\Hrt^{\textup{trunc}}(P)=5.00$ corresponding to Problem~\ref{prob:returntimeentropy} in the SF map.}\label{fig:SFreturn}
\end{figure}

Since the closed form of the return time entropy $\Hrt(P)$ is not available, we truncate the infinite series and obtain a truncated objective function for computation and optimization purposes. In order to control the difference between the original and the truncated objective function, for a given accuracy level $\eta$, we use Markov inequality to select the truncation position $N_\eta$ such that the discarded tail return time probabilities for any location is under $\eta$. Then, we define the truncated objective function as
\begin{equation*}
\Hrt^{\textup{trunc}}(P)=-\sum_{i=1}^n\pi_i\sum\limits_{k=1}^{N_\eta} F_k(i,i)\log F_k(i,i).
\end{equation*}
The objective function in Problem~\ref{prob:returntimeentropy} is replaced with $\Hrt^{\textup{trunc}}(P)$ when solving the optimization problem. Based on the gradient formulas for $\Hrt^{\textup{trunc}}(P)$ \cite[Lemma 18]{XD-MG-FB:17o}, we utilize the gradient projection algorithm to compute a locally optimal solution. We plot the found optimal solution for SF map in~Fig.~\ref{fig:SFreturn}, where $N_\eta$ is taken to be $2292$ so that the discarded probability $\eta=0.1$. We refer the readers to \cite{XD-MG-FB:17o} for a performance comparison with other Markov chains in capturing intruders with certain behaviors.

\section{SUMMARY}\label{sec:conclusion}
In this paper, we presented a stochastic approach to the design of robotic surveillance algorithms over graphs. This approach exhibits numerous promising properties that are relevant in the surveillance applications. We formulated various optimization problems in different settings when particular metrics are concerned. These formulations accommodate scenarios where different characteristics of surveillance strategies are desirable. To facilitate the computation of many of the strategies discussed in this paper, we developed software packages available in both MATLAB and Julia \cite{HW:08-github}.

Stochastic robotic surveillance is still an active area of research and many problems remain open. First, regarding the modeling, it is becoming increasingly popular to incorporate adversary models in the problem formulation. By doing so, we may expect improved
surveillance performance due to the tailored design against the adversarial behaviors. Game theory is a common tool that shall fit in the picture well. Secondly, there is still a lack of efficient algorithms to compute surveillance strategies with optimality guarantees. Although effective heuristics have been proposed in different scenarios, progress in advancing customized computational algorithms is critical. 
Finally, most of the problems we discussed in this paper involve a single surveillance robot. Studies on how to exploit robot teams and design cooperative surveillance systems are of tremendous value in both research and practice.
\section*{DISCLOSURE STATEMENT}
The authors are not aware of any affiliations, memberships, funding, or financial holdings that
might be perceived as affecting the objectivity of this review. 

\section*{ACKNOWLEDGMENTS}
This work has been supported in part by Air Force Office of Scientific
Research award FA9550-15-1-0138. The authors thank
Pushkarini Agharkar,
Andrea Carron,
Mishel George,
Saber Jafarpour, and 
Rushabh Patel
for fruitful past collaborations on the themes of this article.

\bibliographystyle{plainurl+isbn}
\bibliography{alias,Main,FB}

\end{document}